\documentclass[a4paper,11pt,reqno]{amsart}
\usepackage{amssymb,amsmath}

\begin{document}

\title[Iterating lowering operators]{Iterating lowering operators}

\author{Vladimir Shchigolev}

\address{Moscow State University, Russia}

\email{shchigolev\_vladimir@yahoo.com}


\def\K{\mathcal K}
\def\B{\mathfrak B}
\def\C{\mathfrak C}
\def\Z{\mathbb Z}
\def\GL{\mathop{\rm GL}\nolimits}

\def\im{\mathop{\rm Im}}
\def\res{\mathop{\rm res}\nolimits}
\def\Proof{\noindent{\bf Proof. }}
\def\endproof{\hfill$\square$}

\renewcommand{\le}{\leqslant}
\renewcommand{\ge}{\geqslant}
\renewcommand\emptyset{\varnothing}
\renewcommand\epsilon{\varepsilon}
\renewcommand{\labelenumi}{{\rm \theenumi}}
\renewcommand{\theenumi}{{\rm(\roman{enumi})}}

\makeatletter
\def\pmod#1{\allowbreak\if@display\mkern6mu\else\mkern6mu\fi({\operator@font mod}\mkern6mu#1)}
\makeatother

\newtheorem{theorem}{Theorem}
\newtheorem{lemma}[theorem]{Lemma}
\newtheorem{conjecture}[theorem]{Conjecture}
\newtheorem{problem}[theorem]{Problem}
\newtheorem{corollary}[theorem]{Corollary}
\newtheorem{definition}[theorem]{Definition}
\newtheorem{proposition}[theorem]{Proposition}

\newcounter{remark}
\newcommand{\remark}{\par\refstepcounter{remark}%
{\it Remark \arabic{remark}.} }
\renewcommand{\theremark}{\arabic{remark}}

\begin{abstract}
For an algebraically closed base field of positive characteristic,
an algorithm to construct some non-zero $\GL(n-1)$-high weight vectors
of irreducible rational $\GL(n)$-modules is suggested.
It is based on the criterion proved in this paper
for the existence of a set $A$
such that $S_{i,j}(A)f_{\mu,\lambda}$ is
a non-zero $\GL(n-1)$-high weight vector,
where $S_{i,j}(A)$ is Kleshchev's lowering operator
and $f_{\mu,\lambda}$ is a non-zero $\GL(n-1)$-high weight vector
of weight $\mu$ of the costandard $\GL(n)$-module $\nabla_n(\lambda)$
with highest weight $\lambda$.
\end{abstract}

\maketitle

\section{Introduction}

Classical lowering operators were introduced by Carter in~\cite{Carter1}.
Kleshchev used them in~\cite{Kleshchev2} to define generalized lowering
operators. Following~\cite{Brundan_quantum} and~\cite{Kleshchev_gjs11},
we denote these operators by $S_{i,j}(A)$. Kleshchev's lowering operators are
useful in constructing $\GL(n-1)$-high weight vectors from the first level of
irreducible rational $\GL(n)$-modules. In fact, \cite[Theorem~4.2]{Kleshchev2}
shows that every such vector has the form $S_{i,n}(A)v_+$, where $v_+$
is the $\GL(n)$-high weight vector. A natural idea is to continue to apply
lowering operators $S_{i,j}(A)$ to the $\GL(n-1)$-high weight vectors already
obtained in order to construct new $\GL(n-1)$-high weight vectors belonging to
higher levels. For example, this method (for $j=n$) was used
in~\cite{Kleshchev_gjs11} to construct all $\GL(n-1)$-high weight vectors of
irreducible modules $L_n(\lambda)$, where $\lambda$ is a generalized
Jantzen-Seitz weight. The main aim of this paper is to find all $\GL(n-1)$-high
weight vectors that can be constructed in this way (see Theorem~\ref{theorem:3}
and Remark~\ref{remark:1} for removing one node and Theorems~\ref{theorem:5}
and~\ref{theorem:6} for moving one node).

Let $K$ be an algebraically closed field of characteristic $p>0$
and $\GL(m)$ denote the group of invertible $m\times m$-matrices over $K$.
We generally follow the notations of~\cite{Kleshchev_gjs11}
and~\cite{Brundan_quantum} and actually work with hyperalgebras rather than
algebraic groups.
For the connection between representations of the latter two,
we refer the reader to~\cite{Jantzen1}.
Let $U(m,\Z)$ denote the $\Z$-subalgebra of the universal
enveloping algebra $U(m,{\mathbb C})$ of the Lie algebra
${\mathfrak gl}(m,{\mathbb C})$ that is generated by the identity element
and
$$
\begin{array}{l}
X_{i,j}^{(r)}:=\frac{(X_{i,j})^r}{r!}\mbox{ for }1\le i,\; j\le m, i\ne j\mbox{ and }r\ge1;\\[3pt]
\binom{X_{i,i}}{r}:=\frac{X_{i,i}(X_{i,i}-1)\cdots(X_{i,i}-r+1)}{r!}\mbox{ for }1\le i\le m\mbox{ and }r\ge1,
\end{array}
$$
where $X_{i,j}$ denotes the $m\times m$-matrix
with $1$ in the $ij$-entry and zeros elsewhere.
We define the {\it hyperalgebra} $U(m)$ to be
$U(m,\Z)\otimes_\Z K$. For $1\le i<j\le m$ we denote by
$E^{(r)}_{i,j}$ and $F^{(r)}_{i,j}$ the images of
$X^{(r)}_{i,j}$ and $X^{(r)}_{j,i}$ respectively and
for $1\le i\le m$ denote by $\binom{H_i}{r}$
the image of $\binom{X_{i,i}}{r}$
under the above base change.
If $r=1$ then
we omit the superscripts in the above definitions and
write $H_i$ for $\binom{H_i}{1}$.
We also put $E^{(r)}_i:=E^{(r)}_{i,i+1}$ and $F^{(r)}_{i,i}:=1$.

Let $U^0(m)$ denote the subalgebra of $U(m)$ generated by
$1$ and $\binom{H_i}{r}$ for $1\le i\le m$ and $r\ge1$
and $X^+(m)$ denote the set of
integer sequences $(\lambda_1,\ldots,\lambda_m)$ such that $\lambda_1\ge\cdots\ge\lambda_m$.
We say that a vector $v$ of a $U(m)$-module {\it has weight} $\lambda\in X^+(m)$
if $\binom{H_i}{r}v=\binom{\lambda_i}{r}v$ for any $1\le i\le m$ and $r\ge 1$.
If moreover $E^{(r)}_iv=0$ for any $1\le i<m$ and $r\ge 1$, then
we say that $v$ is a {\it $U(m)$-high weight vector}.

Throughout $[i..j]$, $(i..j]$,  $[i..j)$, $(i..j)$
denote the sets
$\{a\in\Z:i\le a\le j\}$,
$\{a\in\Z:i<   a\le j\}$,
$\{a\in\Z:i\le a<   j\}$,
$\{a\in\Z:i<   a<   j\}$ respectively.
For any condition $\mathcal P$, let $\delta_{\mathcal P}$ be $1$
if $\mathcal P$ is true and $0$ if it is false.
Given a pair of integers $(i,j)$, let $\res_p(i,j)$
denote $(i-j)+p\Z$, which is an element of $\Z/p\Z$.
For any set $A\subset\Z$ and two integers $i\le j$,
let $A_{i..j}=\{a\in A:i<a<j\}$. If moreover $A\subset(i..j)$ then
we put $F_{i,j}^A=F_{a_0,a_1}\cdots F_{a_k,a_{k+1}}$, where
$A\cup\{i,j\}=\{a_0<\cdots<a_{k+1}\}$.
Thus $F^\emptyset_{i,j}=F_{i,j}$.
For $i<j$ and $A\subset(i..j)$, the lowering operator $S_{i,j}(A)$ is
defined as (see~\cite[Remark 4.8]{Brundan_quantum})
$$
S_{i,j}(A):=\sum_{B\subset(i..j)}F_{i,j}^BH_{i,j}(A,B).
$$
In this formula, $H_{i,j}(A,B)$ is the element of $U^0(m)$
obtained by evaluating the rational expression
$$
{\mathcal H}_{i,j}(A,B):=\sum_{D\subset B\setminus A}(-1)^{|D|}
\dfrac{\prod\limits_{t\in A}(x_t-x_{D_i(t)})}{\prod\limits_{t\in B}(x_t-x_{D_i(t)})},
$$
where $D_i(t)=\max\{s\in D\cup\{i\}:s<t\}$,
at $x_k:=k-H_k$.
Elements $H_{i,j}(A,B)$ are well defined,
since ${\mathcal H}_{i,j}(A,B)\in\Z[x_i,\ldots,x_{j-1}]$,
which is proved in~\cite[Lemma 4.6(i)]{Brundan_quantum}.
We additionally assume that $S_{i,i}(\varnothing)=1$.

Quite easy proofs of all the properties of
the operators $S_{i,j}(A)$ we need here can be found in~\cite{Brundan_quantum},
where the specialization $v\mapsto1$ should be made.

In this paper, we work with costandard modules $\nabla_n(\lambda)$,
where $\lambda\in X^+(n)$, and its non-zero
$U(n-1)$-high weight vectors $f_{\mu,\lambda}$,
where $\mu\in X^+(n-1)$ and $\lambda_i\ge\mu_i\ge\lambda_{i+1}$ for $1\le i<n$.
If the last conditions hold we write $\mu\longleftarrow\lambda$.
We also denote the element $f_{\bar\lambda,\lambda}$,
where $\bar\lambda=(\lambda_1,\ldots,\lambda_{n-1})$, by $f_\lambda$.
It is a $U(n)$-high weight vector
generating the simple submodule $L_n(\lambda)$ of $\nabla_n(\lambda)$.
The definitions of all these objects can be found in~\cite{Kleshchev_gjs11}.
Moreover using~\cite[Lemma 2.6(ii)]{Kleshchev_gjs11} and
multiplication by a suitable power of
the determinant representation of $\GL(n)$,
we may assume that
$f_\lambda$ and $f_{\mu,\lambda}$,
where $\mu\longleftarrow\lambda$ and $a_i:=\sum_{s=1}^i(\lambda_s-\mu_s)$,
are chosen so that
$E_1^{(a_1)}\cdots E_{n-1}^{(a_{n-1})}f_{\mu,\lambda}=f_\lambda$.

\section{Graph of sequences}
For the remainder of this paper, we fix an integer $n>1$ and
weights $\lambda\in X^+(n)$, $\mu\in X^+(n-1)$ such that $\mu\longleftarrow\lambda$.
For $i=1$, \ldots, $n-1$, we put $a_i:=\sum_{j=1}^i(\lambda_j-\mu_j)$.
The following formulas can easily be checked by calculations in $U(n,\Z)$.

\begin{lemma}\label{lemma:1}
Let $1\le i<j\le n$, $1\le l<n$, $m\ge1$ and $A\subset(i..j)$. We have
\begin{enumerate}
\item $E_l^{(m)}F_{i,j}^A=F_{i,j}^AE_l^{(m)}$
      if $l\notin A\cup\{i\}$ and $l+1\notin A\cup\{j\}$;
\item $E_l^{(m)}F_{i,j}^A=F_{i,j}^AE_l^{(m)}-F_{i,l}^{A_{i..l}}F_{l+1,j}^{A_{l+1..j}}E_l^{(m-1)}$
      if $l\in A\cup\{i\}$ and $l+1\notin A\cup\{j\}$;
\item $E_l^{(m)}F_{i,j}^A=F_{i,j}^AE_l^{(m)}+F_{i,l}^{A_{i..l}}F_{l+1,j}^{A_{l+1..j}}E_l^{(m-1)}$
      if $l\notin A\cup\{i\}$ and $l+1\in A\cup\{j\}$;
\item $E_l^{(m)}F_{i,j}^A=F_{i,j}^AE_l^{(m)}+
      F_{i,l}^{A_{i..l}}(H_l-H_{l+1}+1-m)F_{l+1,j}^{A_{l+1..j}}E_l^{(m-1)}$
      if $l\in A\cup\{i\}$ and $l+1\in A\cup\{j\}$.
\end{enumerate}
\end{lemma}
We shall use the abbreviation $E(i,j)=E_i\cdots E_j$.
Let $1\le i\le k\le j\le n$ and $A\subset(i..j)$.
It follows from Lemma~\ref{lemma:1} that
$E(k,j-1)S_{i,j}(A)=u_kE_k+\cdots+u_{j-1}E_{j-1}+M^k_{i,j}(A)$,
where $u_k$, \ldots, $u_{j-1}\in U(n)$ and
$M^k_{i,j}(A)$ is a linear combination
of elements of the form $F_{i,k}^BH$, where $H\in U^0(n)$.
In what follows, we stipulate that any not necessarily commutative
product of the form $\prod_{i\in A}x_i$,
where $A=\{a_1<\cdots<a_m\}\subset\Z$, equals $x_{a_1}\cdots x_{a_m}$.

\begin{lemma}\label{lemma:2}
Given integers $1\le i_1<j_1<\cdots<i_{s-1}<j_{s-1}<i_s<j_s\le n$,
sets $A_1\subset(i_1..j_1)$, \ldots, $A_s\subset(i_s..j_s)$
and integers $k_1,\ldots,k_s$ such that $i_t\le k_t\le j_t$ for $t=1$, \ldots, $s$
and $j_s=n$ implies $k_s=n$, we put
$$
v=E(k_1,j_1-1)S_{i_1,j_1}(A_1)\cdots E(k_s,j_s-1)S_{i_s,j_s}(A_s)f_{\mu,\lambda}.
$$
Then we have
\begin{enumerate}
\item\label{lemma:2:p:1}
$v=X_1\cdots X_sf_{\mu,\lambda}$, where
each $X_t$ is either $E(k_t,j_t-1)S_{i_t,j_t}(A_t)$ or
$M^{k_t}_{i_t,j_t}(A_t)$;
\item\label{lemma:2:p:2} $E_l^{(m)}v=0$ if $1\le l<n-1$ and $m\ge2$;
\item\label{lemma:2:p:3} $E_l^{(m)}v=0$ if $m\ge1$ and
     $l\in[1..n-1)\setminus\bigl([i_1..k_1)\cup\cdots\cup[i_s..k_s)\bigl)$;
\item\label{lemma:2:p:4} If $i_t<k_t<n$ then
$$
\begin{array}{l}
\displaystyle E_{k_t-1}v=
\left(\prod_{r=1}^{t-1}E(k_r,j_r-1)S_{i_r,j_r}(A_r)\right)
E(k_t-1,j_t-1)S_{i_t,j_t}(A_t)\\
\displaystyle\times\left(\prod_{r=t+1}^sE(k_r,j_r-1)S_{i_r,j_r}(A_r)\right)f_{\mu,\lambda};
\end{array}
$$
\item\label{lemma:2:p:5} If $l\in[i_t..k_t-1)$ then
$$
\begin{array}{l}
\displaystyle E_lv=
c\left(\prod_{r=1}^{t-1}E(k_r,j_r-1)S_{i_r,j_r}(A_r)\right)
S_{i_t,l}((A_t)_{i_t..l})\\
\displaystyle\times E(k_t,j_t-1)S_{l+1,j_t}((A_t)_{l+1..j_t})
\left(\prod_{r=t+1}^sE(k_r,j_r-1)S_{i_r,j_r}(A_r)\right)f_{\mu,\lambda},
\end{array}
$$
where $c=0$ except the case
$l\in A_t\cup\{i_t\}$, $l+1\notin A_t$,
in which $c=-1$.
\end{enumerate}
\end{lemma}
\Proof~\ref{lemma:2:p:1}
Applying Lemma~\ref{lemma:1}, we prove by induction on $t$
(starting from $t=s$) that
$$
\begin{array}{l}
v=E(k_1,j_1-1)S_{i_1,j_1}(A_1)\cdots
E(k_{t-1},j_{t-1}-1)S_{i_{t-1},j_{t-1}}(A_{t-1})\\[3pt]
\times M^{k_t}_{i_t,j_t}(A_t)\cdots M^{k_s}_{i_s,j_s}(A_s)f_{\mu,\lambda}.
\end{array}
$$
Using this formula for $t=1$,
we obtain the required result by induction on $s$.

\ref{lemma:2:p:2}, \ref{lemma:2:p:3}
follow from part~\ref{lemma:2:p:1} for $X_t=M^{k_t}_{i_t,j_t}(A_t)$
and Lemma~\ref{lemma:1}.

\ref{lemma:2:p:4}
Applying part~\ref{lemma:2:p:1} (possibly for different parameters),
we get
$$
\begin{array}{l}
E_{k_t-1}v=
E_{k_t-1}M^{k_1}_{i_1,j_1}(A_1)\cdots M^{k_{t-1}}_{i_{t-1},j_{t-1}}(A_{t-1})\\[3pt]
\times E(k_t,j_t{-}1)S_{i_t,j_t}(A_t)\cdots E(k_s,j_s{-}1)S_{i_s,j_s}(A_s)f_{\mu,\lambda}\\[3pt]
=M^{k_1}_{i_1,j_1}(A_1)\cdots M^{k_{t-1}}_{i_{t-1},j_{t-1}}(A_{t-1})
E(k_t-1,j_t-1)S_{i_t,j_t}(A_t)\\[3pt]
\times E(k_{t+1},j_{t+1}-1)S_{i_{t+1},j_{t+1}}(A_{t+1})
\cdots E(k_s,j_s-1)S_{i_s,j_s}(A_s)f_{\mu,\lambda}.
\end{array}
$$
Now the required formula follows from part~\ref{lemma:2:p:1}.

\ref{lemma:2:p:5}
Since $E_l$ and $E(k_t,j_t-1)$ commute in this case,
we get by~\cite[4.11(i),(ii)]{Brundan_quantum} and
parts~\ref{lemma:2:p:1},\ref{lemma:2:p:2} of the current lemma that
$$
\begin{array}{l}
E_lv=M^{k_1}_{i_1,j_1}(A_1)\cdots M^{k_{t-1}}_{i_{t-1},j_{t-1}}(A_{t-1})
E(k_t,j_t-1)E_lS_{i_t,j_t}(A_t)\\[3pt]
\times E(k_{t+1},j_{t+1}-1)S_{i_{t+1},j_{t+1}}(A_{t+1})
\cdots E(k_s,j_s-1)S_{i_s,j_s}(A_s)f_{\mu,\lambda}=\\[3pt]
cM^{k_1}_{i_1,j_1}(A_1)\cdots M^{k_{t-1}}_{i_{t-1},j_{t-1}}(A_{t-1})
S_{i_t,l}((A_t)_{i_t..l})E(k_t,j_t-1)S_{l+1,j_t}((A_t)_{l+1..j_t})\\[3pt]
\times E(k_{t+1},j_{t+1}-1)S_{i_{t+1},j_{t+1}}(A_{t+1})
\cdots E(k_s,j_s-1)S_{i_s,j_s}(A_s)f_{\mu,\lambda}.
\end{array}
\!
\!
\!
$$
Now the required formula follows similarly to~\ref{lemma:2:p:4}.
\endproof

For $1\le i<j\le n$ and $A\subset(i..j)$,
we define the polynomial $\K_{i,j}(A)$ of
$\Z[x_i,\ldots,x_{j-1},y_{i+1},\ldots,y_j]$
as in~\cite[4.12]{Brundan_quantum} by the formula
$$
\K_{i,j}(A):=\sum_{B\subset(i..j)}\left({\mathcal H}_{i,j}(A,B)\prod_{t\in B\cup\{i\}}(y_{t+1}-x_t)\right).
$$
We define
$H_{i,j}^\mu(A,B)$ by evaluating
${\mathcal H}_{i,j}(A,B)$ at $x_q:=\res_p(q,\mu_q)$
and define
$K_{i,j}^{\mu,\lambda,k}(A)$ by evaluating $\K_{i,j}(A)$ at
\begin{equation}\label{equation:0.5}
\begin{array}{ll}
x_q:=\res_p(q,\mu_q)&\mbox{for }1\le q<n,\\
y_q:=\res_p(q,\lambda_q+1)&\mbox{for }1<q\le k,\\
y_q:=\res_p(q,\mu_q+1)&\mbox{for }k<q<n,
\end{array}
\end{equation}
where $1+\delta_{j=n}(n-1)\le k\le n$.
For $1\le i\le t<n$ and
$1+\delta_{t+1=n}(n-1)\le k\le n$, let $B^{\mu,\lambda,k}(i,t)$ denote
the element of $\Z/p\Z$ obtained from $y_{t+1}-x_i$
by substitution~(\ref{equation:0.5}).
We also abbreviate
$K_{i,j}^{\mu,\lambda}(A):=K_{i,j}^{\mu,\lambda,n}(A)$ and
$B^{\mu,\lambda}(i,t):=B^{\mu,\lambda,n}(i,t)$.

\remark\label{remark:0.5}
Clearly
$B^{\mu,\lambda,k}(i,t)=t-i+\mu_i-\mu_{t+1}$ for $k\le t$ and
$B^{\mu,\lambda,k}(i,t)=t-i+\mu_i-\lambda_{t+1}$ for $k>t$.
In particular, $B^{\mu,\lambda,k}(i,t)=B^{\mu,\lambda,i}(i,t)$ for $k\le i$ and
$B^{\mu,\lambda,k}(i,t)=B^{\mu,\lambda,t+1}(i,t)$ for $k>t$.

The next result is actually proved in~\cite[Proposition 4.5]{Kleshchev_gjs11}.
Recall that we have defined $a_t=\sum_{j=1}^t(\lambda_j-\mu_j)$.
\begin{proposition}\label{proposition:1}
Given integers $1\le d_1<d'_1\le d_2<d'_2\le\cdots \le d_r<d'_r\le n$,
we have
$$
\left(\prod_{t=1}^{n-1}E_t^{(a_t+\delta_{t\in G})}\right)
F_{d_1,d'_1}\cdots F_{d_r,d'_r}f_{\mu,\lambda}=
\prod_{q=1}^r(\mu_{d_q}-\lambda_{d_q+1})f_\lambda,
$$
where $G=[d_1..d'_1)\cup\cdots\cup[d_r..d'_r)$.
\end{proposition}

\begin{lemma}\label{lemma:3}
Under the hypothesis of Lemma~\ref{lemma:2}, we have
$$
\left(\prod_{t=1}^{n-1}E_t^{(a_t+\delta_{t\in G})}\right)v=
K_{i_1,j_1}^{\mu,\lambda,k_1}(A_1)\cdots K_{i_s,j_s}^{\mu,\lambda,k_s}(A_s)f_\lambda,
$$
where $G=[i_1..k_1)\cup\cdots\cup[i_s..k_s)$.
\end{lemma}
\Proof
By Lemma~\ref{lemma:1}, we have
$E(k_t,j_t-1)F_{i_t,j_t}^B\equiv F_{i_t,k_t}^{B_{i_t..k_t}}\prod_{q\in B\cup\{i_t\},q\ge k_t}(H_q-H_{q+1})$
modulo the left ideal of $U(n)$ generated by
$E_{k_t},\ldots,E_{j_t-1}$.
Thus taking into account~\cite[Remark 4.8]{Brundan_quantum}, we get
\begin{equation}\label{equation:1}
v=\prod_{t=1}^s\sum_{B_t\subset(i_t..j_t)}\left(H_{i_t,j_t}^\mu(A_t,B_t)F_{i_t,k_t}^{(B_t)_{i_t..k_t}}%
\prod_{\genfrac{}{}{0pt}{}{q\in B_t\cup\{i_t\}}{q\ge k_t}}(\mu_q-\mu_{q+1})\right)f_{\mu,\lambda}.
\end{equation}
By Proposition~\ref{proposition:1}, we have
$$
\left(\prod_{t=1}^{n-1}E_t^{(a_t+\delta_{t\in G})}\right)
F_{i_1,k_1}^{(B_1)_{i_1..k_1}}\cdots F_{i_s,k_s}^{(B_s)_{i_s..k_s}}f_{\mu,\lambda}=
\prod_{t=1}^s\prod_{\genfrac{}{}{0pt}{}{q\in B_t\cup\{i_t\}}{q<k_t}}(\mu_q-\lambda_{q+1})f_\lambda.
$$
Substituting this into~(\ref{equation:1}) completes the proof.
\endproof

Let $V_n$ be the set of all sequences
$x=\bigl((i_1,k_1,j_1,A_1)$, \ldots, $(i_s,k_s,j_s,A_s)\bigl)$
such that
$$
\begin{array}{ll}
1\le i_1<j_1<\cdots<i_s<j_s\le n;&
A_1\subset(i_1..j_1),\ldots,A_s\subset(i_s..j_s);\\[3pt]
i_1\le k_1\le j_1,\ldots,i_s\le k_s\le j_s;& j_s=n\mbox{ implies }k_s=n.
\end{array}
$$
Moreover, we put $\Phi(x):=E(k_1,j_1{-}1)S_{i_1,j_1}(A_1)\cdots E(k_s,j_s{-}1)S_{i_s,j_s}(A_s)$
and $K^{\mu,\lambda}(x):=K_{i_1,j_1}^{\mu,\lambda,k_1}(A_1)\cdots K_{i_s,j_s}^{\mu,\lambda,k_s}(A_s)$.
In what follows, we assume that the product of two finite sequences
$a=(a_1,\ldots,a_s)$ and $b=(b_1,\ldots,b_t)$ equals
$ab=(a_1,\ldots,a_s,b_1,\ldots,b_t)$.

Let $x,x'\in V_n$. We write $x\stackrel{l}{\longrightarrow}x'$ if
there exists a representation $x=a\bigl((i,k,j,A)\bigl)b$ such that
one of the following conditions holds:
\begin{itemize}
\item $x'=a\bigl((i,k-1,j,A)\bigl)b$, $l=k-1$, $i<k<n$;
\item $x'=a\bigl((i+1,k,j,A)\bigl)b$, $l=i$, $i+1\notin A$, $i<k-1$;
\item $x'=a\bigl((i,l,l,A_{i..l}),(l+1,k,j,A_{l+1..j})\bigl)b$,
      $l\in(i..k-1)$, $l\in A$, $l+1\notin A$.
\end{itemize}

The above definitions are made exactly to ensure the following property.
\begin{lemma}\label{lemma:3.5}
Let $x,x'\in V_n$. If $x\stackrel{l}{\longrightarrow}x'$ then
$E_l\Phi(x)f_{\mu,\lambda}=\pm\Phi(x')f_{\mu,\lambda}$.
\end{lemma}
\noindent
{\bf Proof} follows directly from
Lemma~\ref{lemma:2}\ref{lemma:2:p:4},\ref{lemma:2:p:5}.\endproof

We say that $x'$ {\it follows} from $x$ if
there are $x_0,\ldots,x_m\in V_n$ and integers $l_0,\ldots,l_{m-1}$
such that
$x=x_0$, $x'=x_m$ and $x_t\stackrel{l_t}{\longrightarrow}x_{t+1}$
for $0\le t<m$.
In particular, every element of $V_n$ follows from itself.

\begin{theorem}\label{theorem:1}
Let $x\in V_n$.
The equality $\Phi(x)f_{\mu,\lambda}=0$ holds if and only if
$K^{\mu,\lambda}(x')=0$ for any $x'$ following from $x$.
\end{theorem}
\Proof
It follows from Lemmas~\ref{lemma:3.5} and~\ref{lemma:3} that
$\Phi(x)f_{\mu,\lambda}=0$ implies $K^{\mu,\lambda}(x')=0$
for any $x'$ following from $x$.

Let $x=\bigl((i_1,k_1,j_1,A_1),\ldots,(i_s,k_s,j_s,A_s)\bigl)$.
We prove the reverse implication by induction on $\sum_{t=1}^s(k_t-i_t)$.
The induction starts by noting that this sum is always non-negative.
So we suppose that the reverse implication is true for smaller
values of this sum.
By Lemma~\ref{lemma:2}\ref{lemma:2:p:2},\ref{lemma:2:p:3},
we get $E_l^{(m)}\Phi(x)f_{\mu,\lambda}=0$ if $l<n-1$ and $m>1$ or if
$m\ge1$ and
$l\in[1..n-1)\setminus\bigl([i_1..k_1)\cup\cdots\cup[i_s..k_s)\bigl)$.

However $E_l\Phi(x)f_{\mu,\lambda}=0$ also for $l\in[1..n-1)\cap\bigl([i_1..k_1)\cup\cdots\cup[i_s..k_s)\bigl)$
by Lemma~\ref{lemma:3.5} and the inductive hypothesis.
Thus $\Phi(x)f_{\mu,\lambda}$ is a $U(n-1)$-high weight vector of
weight $\nu=\mu-\sum_{t=1}^s(\varepsilon_{i_t}-\varepsilon_{k_t})$,
where $\varepsilon_i=(0^{i-1},1,0^{n-1-i})$ for $i<n$ and $\varepsilon_n=(0^{n-1})$.
It follows from~\cite[Corollary 3.3]{Kleshchev_gjs11} that
$\Phi(x)f_{\mu,\lambda}=0$ if $\nu\longleftarrow\lambda$ does not hold
and that $\Phi(x)f_{\mu,\lambda}=cf_{\nu,\lambda}$ for some $c\in K$ if
$\nu\longleftarrow\lambda$. We need to consider only the latter case.
By the last equation of the introduction
and Lemma~\ref{lemma:3},
we have
$cf_\lambda=X(cf_{\nu,\lambda})=X\Phi(x)f_{\mu,\lambda}=K^{\mu,\lambda}(x)f_\lambda=0$,
where $X=\prod_{t=1}^{n-1}E_t^{(a_t+\delta_{t\in G})}$ and
$G=[i_1..k_1)\cup\cdots\cup[i_s..k_s)$.
Hence $c=0$ and $\Phi(x)f_{\mu,\lambda}=0$. \endproof

The next corollary follows from Theorem~\ref{theorem:1} and
the following simple fact:
if $x\in V_n$ and $x=x_1x_2$ then $x'$ follows from $x$ if and only if
there are sequences $x'_1$ and $x'_2$ following from
$x_1$ and $x_2$ respectively such that $x'=x'_1x'_2$.

\begin{corollary}\label{corollary:1}
Let $x\in V_n$ and $x=x_1x_2$.
Then $\Phi(x)f_{\mu,\lambda}=0$ if and only if
$\Phi(x_1)f_{\mu,\lambda}=0$ or $\Phi(x_2)f_{\mu,\lambda}=0$.
\end{corollary}

\section{Removing one node}

We say that a map $\theta:A\to\Z$, where $A\subset\Z$,
is {\it weakly increasing} ({\it weakly decreasing}) if
$\theta(a)\ge a$ (resp. $\theta(a)\le a$) for any $a\in A$.
We need the following facts about the polynomials $\K_{i,j}(A)$.

\begin{proposition}\label{proposition:2}
Let $1\le i<j\le n$, $1+\delta_{j=n}(n-1)\le k\le n$, $A\subset(i..j)$ and
there exists
a weakly increasing injection $\theta:(i..j)\setminus A\to(i..j)$
such that $B^{\mu,\lambda,k}(t,\theta(t))=0$ for any $t\in(i..j)\setminus A$.
Then $$K_{i,j}^{\mu,\lambda,k}(A)=\prod_{t\in[i..j)\setminus\im\theta}B^{\mu,\lambda,k}(i,t).$$
\end{proposition}
\noindent
\Proof The result is obtained from~\cite[Lemma 4.4]{Kleshchev_gjs11}
by substitution~(\ref{equation:0.5}). \endproof

\begin{lemma}\label{lemma:4}
For $i<j-1$ and $A\subset(i..j)$, we have
\begin{enumerate}
\item\label{lemma:4:p:1} $\K_{i,j}(A)=\K_{i,j-1}(A)$ if $j-1\notin A$;
\item\label{lemma:4:p:2} $\K_{i,j}(A)=\K_{i,j-1}(A\setminus\{j-1\})(y_j-x_k)+
                         \delta_{k\ne i}\K_{i,j-1}(\{k\}\cup A\setminus\{j-1\})$,
                         where $k=\max[i..j)\setminus A$, if $j-1\in A$.
\end{enumerate}
\end{lemma}
\Proof
We put $\bar A=(i..j)\setminus A$.
In this proof, we use~\cite[Lemma 4.13(i)]{Brundan_quantum} for
a self-contained form of $\K_{i,j}(A)$ and
the following notation of~\cite{Brundan_quantum}:
if $D\subset(i..j)$ and $k>i$ then $D_i(k)=\max\{t\in D\cup\{i\}:t<k\}$.

\ref{lemma:4:p:1}
If $D\subset\bar A\setminus\{j-1\}$ then
$(D\cup\{j-1\})_i(t)=D_i(t)$ for $t<j$,
$(D\cup\{j-1\})_i(j)=j-1$ and $D_i(j)=D_i(j-1)$.
Hence we get
$$
\begin{array}{l}
\K_{i,j}(A)=
{\displaystyle\sum_{D\subset\bar A\setminus\{j-1\}}(-1)^{|D|}}
\left(
\frac{\prod\limits_{t\in(i..j]}(y_t-x_{D_i(t)})}{\prod\limits_{t\in\bar A}(x_t-x_{D_i(t)})}-
\frac{\prod\limits_{t\in(i..j]}(y_t-x_{(D\cup\{j-1\})_i(t)})}{\prod\limits_{t\in\bar A}(x_t-x_{(D\cup\{j-1\})_i(t)})}
\right)=\\[20pt]
{\displaystyle\sum_{D\subset\bar A\setminus\{j-1\}}(-1)^{|D|}}
\left(
\frac{\prod\limits_{t\in(i..j-1]}(y_t-x_{D_i(t)})}{\prod\limits_{t\in\bar A\setminus\{j-1\}}(x_t-x_{D_i(t)})}
\frac{(y_j-x_{D_i(j-1)})-(y_j-x_{j-1})}{x_{j-1}-x_{D_i(j-1)}}
\right)
=\K_{i,j-1}(A).
\end{array}
\!\!\!\!
$$

\ref{lemma:4:p:2}
If $k=i$ then $A=(i..j)$,
$\K_{i,j}(A)=\prod_{t\in(i..j]}(y_t-x_i)$,
$\K_{i,j-1}(A\setminus\{j-1\})=\prod_{t\in(i..j-1]}(y_t-x_i)$
(by part~\ref{lemma:4:p:1}) and the required formula follows.

Therefore, we consider the case $k\ne i$.
We have
$$
\begin{array}{l}
\K_{i,j}(A)=(y_j-x_k)
{\displaystyle \sum_{D\subset\bar A}(-1)^{|D|}}
\dfrac{\prod_{t\in(i..j-1]}(y_t-x_{D_i(t)})}{\prod_{t\in\bar A}(x_t-x_{D_i(t)})}\\[12pt]
+
{\displaystyle\sum_{D\subset\bar A}(-1)^{|D|}}
(x_k-x_{D_i(j)})
\dfrac{\prod_{t\in(i..j-1]}(y_t-x_{D_i(t)})}{\prod_{t\in\bar A}(x_t-x_{D_i(t)})}.
\end{array}
$$
Part~\ref{lemma:4:p:1} shows that the first sum equals
$\K_{i,j}(A\setminus\{j-1\})$.
Let us look at the second sum. If $k\in D$ then
$D_i(j)=k$ and the summands corresponding to such sets $D$ can be omitted.
If $k\notin D$ then $D_i(j)=D_i(k)$ and this summand equals
$$
(-1)^{|D|}\dfrac{\prod_{t\in(i..j-1]}(y_t-x_{D_i(t)})}
{\prod_{t\in\bar A\setminus\{k\}}(x_t-x_{D_i(t)})}.
$$
Thus the second sum equals $\K_{i,j-1}(\{k\}\cup A\setminus\{j-1\})$. \endproof

Next, we are going to prove the result similar
to~\cite[Proposition 3.2]{Kleshchev2}, where
we replace the $U(n)$-high weight vector $v_+$ by
the $U(n-1)$-high weight vector $f_{\mu,\lambda}$.
The general scheme of proof is borrowed from~\cite[Proposition 3.2]{Kleshchev2},
although some changes are necessary. We shall use Theorem~\ref{theorem:1}
and Lemma~\ref{lemma:4} to make them.
In what follows, we say that a formula
$M=[b_1..c_1]\cup\cdots\cup[b_N..c_N]$ is
{\it the decomposition of $M$
into the union of connected components}
if $b_i\le c_i$ for $1\le i\le N$ and
$c_i<b_{i+1}-1$ for $1\le i<N$.

\begin{definition}\label{definition:1}
Let $1\le i<j\le n$, $M\subset(i..j)$ and
$M=[b_1..c_1]\cup\cdots\cup[b_N..c_N]$ be the decomposition of $M$
into the union of connected components.
We say that $M$ satisfies the condition $\pi_{i,j}^{\mu,\lambda}(v)$ if
$1\le v\le N+1$ and for any $k=1+\delta_{b_v-1=n}(n-1),\ldots,n$
there exists
a weakly increasing injection
$\theta_k:\{i\}\cup[b_1..c_1]\cup\cdots\cup[b_{v-1}..c_{v-1}]\to[i..b_v-1)$
such that $B^{\mu,\lambda,k}(x,\theta_k(x))=0$ for any admissible $x$,
where we assume $b_{N+1}=j+1$.
\end{definition}

\begin{lemma}\label{lemma:5}
Let $1\le i<j\le n$ and $A\subset(i..j)$ be such that
$(i..j)\setminus A$ satisfies $\pi_{i,j}^{\mu,\lambda}(v)$ for some $v$.
Then $K_{i,j}^{\mu,\lambda,k}(A)=0$ for $1+\delta_{j=n}(n-1)\le k\le n$.
\end{lemma}
\Proof
Let $(i..j)\setminus A=[b_1..c_1]\cup\cdots\cup[b_N..c_N]$ be
the decomposition into the union of connected components.
Note that if $v=N+1$, then the required equalities immediately follow
from Proposition~\ref{proposition:2}.

Indeed, take any $k=1+\delta_{j=n}(n-1),\ldots,n$.
Since in this case $b_v-1=j$, Definition~\ref{definition:1}
ensures that there exists
a weakly increasing injection
$\theta_k:\{i\}\cup\left((i..j)\setminus A\right)\to[i..j)$
such that $B^{\mu,\lambda,k}(x,\theta_k(x))=0$ for any admissible $x$.
Taking the restriction of $\theta_k$ to $(i..j)\setminus A$ for
$\theta$ in Proposition~\ref{proposition:2}, we obtain
$$
K_{i,j}^{\mu,\lambda,k}(A)=\prod_{t\in[i..j)\setminus\im\theta}B^{\mu,\lambda,k}(i,t).
$$
The last product equals zero, since $B^{\mu,\lambda,k}(i,\theta_k(i))=0$ and
$\theta_k(i)\in[i..j)\setminus\im\theta$.

Let us prove the lemma by induction on $j-i$.
The case $j-i=1$ follows from the above remark.
Now let $v\le N$, $j-i>1$ and suppose that the lemma
is true for smaller values of this difference.
Take any $k=1+\delta_{j=n}(n-1),\ldots,n$.
By Lemma~\ref{lemma:4}, we have
$$
K_{i,j}^{\mu,\lambda,k}(A)=K_{i,j-1}^{\mu,\lambda,k}(A\setminus\{j-1\})B+
K_{i,j-1}^{\mu,\lambda,k}(\{c_N\}\cup A\setminus\{j-1\})
$$
if $c_N<j-1$ and
$$
K_{i,j}^{\mu,\lambda,k}(A)=K_{i,j-1}^{\mu,\lambda,k}(A)
$$
if $c_N=j-1$, where
$B$ is the element of $\Z/p\Z$ obtained from $y_j-x_{c_N}$
by substitution~(\ref{equation:0.5}).
Clearly, the sets $(i..j-1)\setminus(A\setminus\{j-1\})$ and
$(i..j-1)\setminus(\{c_N\}\cup A\setminus\{j-1\})$ in the former case
and the set
$(i..j-1)\setminus A$
in the latter case satisfy the condition $\pi_{i,j-1}^{\mu,\lambda}(v)$.
\endproof

\begin{theorem}\label{theorem:2}
Let $1\le i<j\le n$ and $A\subset(i..j)$.
Then $S_{i,j}(A)f_{\mu,\lambda}=0$ if and only if
$(i..j)\setminus A$ satisfies $\pi_{i,j}^{\mu,\lambda}(v)$ for some $v$.
\end{theorem}
\Proof
Let $\bar A=(i..j)\setminus A$ and
$\bar A=[b_1..c_1]\cup\cdots\cup[b_N..c_N]$ be
the decomposition into the union of connected components.
We put $x_k=\bigl((i,k,j,A)\bigl)$ for brevity.
It should be kept in mind that $\Phi(x_j)=S_{i,j}(A)$.

We prove the theorem by induction on $|\bar A|$.
Suppose $\bar A=\emptyset$.
Then all the sequences following from $x_j$ are
$x_k$, where $i+\delta_{j=n}(j-i)\le k\le j$.
By Theorem~\ref{theorem:1},
$\Phi(x_j)f_{\mu,\lambda}=0$ if and only if
$K^{\mu,\lambda}(x_k)=0$ for any $k=i+\delta_{j=n}(j-i),\ldots,j$.
Applying Proposition~\ref{proposition:2}, we see that
$\Phi(x_j)f_{\mu,\lambda}=0$ if and only if
for any $k=i+\delta_{j=n}(j-i),\ldots,j$ there is
$t_k\in[i..j)$ such that $B^{\mu,\lambda,k}(i,t_k)=0$.
In view of Remark~\ref{remark:0.5}, this assertion is equivalent to
$\pi_{i,j}^{\mu,\lambda}(1)$.

Now suppose that $\bar A\ne\emptyset$ and that the theorem holds for
smaller values of $|\bar A|$.

{\it``If part''}.
By~\cite[4.11(ii)]{Brundan_quantum} for any $m=1,\ldots,N$,
we have
$E_{b_m-1}S_{i,j}(A)f_{\mu,\lambda}=
-S_{i,b_m-1}(A_{i..b_m-1})$ $S_{b_m,j}(A_{b_m..j})f_{\mu,\lambda}$.
Note that
\begin{equation}\label{equation:1.25}
\begin{array}{l}
A_{i..b_m-1}=(i..b_m-1)\setminus\bigl([b_1..c_1]\cup\cdots\cup[b_{m-1}..c_{m-1}]\bigl),\\[3pt]
A_{b_m..j}=(b_m..j)\setminus\bigl((b_m..c_m]\cup\cdots\cup[b_N..c_N]\bigl).
\end{array}
\end{equation}
If $m\le v-1$ then $(b_m..c_m]\cup\cdots\cup[b_N..c_N]$
satisfies
$\pi_{b_m,j}^{\mu,\lambda}(v-m+1-\delta_{b_m=c_m})$,
whence by the inductive hypothesis $S_{b_m,j}(A_{b_m..j})f_{\mu,\lambda}=0$.
If $m\ge v$ then $i<b_m-1$ and
$[b_1..c_1]\cup\cdots\cup[b_{m-1}..c_{m-1}]$
satisfies
$\pi_{i,b_m-1}^{\mu,\lambda}(v)$,
whence by the inductive hypothesis $S_{i,b_m-1}(A_{i..b_m-1})f_{\mu,\lambda}=0$.
Since the elements $S_{i,b_m-1}(A_{i..b_m-1})$ and $S_{b_m,j}(A_{b_m..j})$
commute, we have in both cases
\begin{equation}\label{equation:1.5}
E_{b_m-1}S_{i,j}(A)f_{\mu,\lambda}=0.
\end{equation}

Let us prove by induction on $s=0,\ldots,j-i$ that in the case $j<n$ the conditions
\begin{equation}\label{equation:2}
K_{i,j}^{\mu,\lambda,j}(A)=0,\quad \ldots,\quad K_{i,j}^{\mu,\lambda,j-s+1}(A)=0,\quad
\Phi(x_{j-s})f_{\mu,\lambda}=0
\end{equation}
imply $\Phi(x_j)f_{\mu,\lambda}=0$.
It is obviously true for $s=0$.
Suppose that $0<s\le j-i$,
conditions~(\ref{equation:2}) hold and
the assertion is true for smaller values of $s$.
By the inductive hypothesis it suffices to prove that
$\Phi(x_{j-s+1})f_{\mu,\lambda}=0$.
Let $x_{j-s+1}\stackrel{l}{\longrightarrow}x'$.
We have either $x'=x_{j-s}$ or $l=b_m-1<j-s$.
Since in the former case $\Phi(x')f_{\mu,\lambda}=0$ by~(\ref{equation:2}),
we shall consider the latter case.
We have
$$
\begin{array}{l}
\Phi(x')f_{\mu,\lambda}=E_{b_m-1}\Phi(x_{j-s+1})f_{\mu,\lambda}=
E_{b_m-1}E(j-s+1,j-1)S_{i,j}(A)f_{\mu,\lambda}\\[6pt]
=E(j-s+1,j-1)E_{b_m-1}S_{i,j}(A)f_{\mu,\lambda}=0.
\end{array}
$$
To obtain the last equality, we used~(\ref{equation:1.5}).
Since $K^{\mu,\lambda}(x_{j-s+1})=K_{i,j}^{\mu,\lambda,j-s+1}(A)=0$,
we get $\Phi(x_{j-s+1})f_{\mu,\lambda}=0$ by Theorem~\ref{theorem:1}.

Note that nothing follows from $x_i$ except itself.
Therefore, applying the above assertion for $s=j-i$ and Theorem~\ref{theorem:1},
we see that to prove $\Phi(x_j)f_{\mu,\lambda}=0$ in the case $j<n$,
it suffices to prove $K_{i,j}^{\mu,\lambda,k}(A)=0$ for $i\le k\le j$.
The last equalities follow from Lemma~\ref{lemma:5}.

If $j=n$ then $x_j\stackrel{l}{\longrightarrow}x'$
holds if and only if $l=b_m-1$, where $1\le m\le N$.
In that case $\Phi(x')f_{\mu,\lambda}=0$ by~(\ref{equation:1.5}).
Therefore, applying Theorem~\ref{theorem:1}, we see that to prove
$\Phi(x_j)f_{\mu,\lambda}=0$ in the case $j=n$,
it suffices to prove
$K_{i,j}^{\mu,\lambda}(A)=0$.
The last equality follows from Lemma~\ref{lemma:5}.

{\it``Only if part''}.
Suppose $\bar A$ satisfies the condition
$\pi_{i,j}^{\mu,\lambda}(v)$ for no $v$.
Multiplying the equality $\Phi(x_j)f_{\mu,\lambda}=0$ by $E_{b_m-1}$,
where $1\le m\le N$, we get
$S_{i,b_m-1}(A_{i..b_m-1})S_{b_m,j}(A_{b_m..j})f_{\mu,\lambda}=0$
according to~\cite[4.11(ii)]{Brundan_quantum}.
By Corollary~\ref{corollary:1}, either
$S_{i,b_m-1}(A_{i..b_m-1})f_{\mu,\lambda}=0$ or
$S_{b_m,j}(A_{b_m..j})f_{\mu,\lambda}=0$.
The former case is impossible since the inductive hypothesis would yield
that $(i..b_m-1)\setminus A_{i..b_m-1}$ satisfies
$\pi_{i,b_m-1}^{\mu,\lambda}(v)$ for some $v\le m$
(see~(\ref{equation:1.25})).
But then
$\bar A$ would satisfy
$\pi_{i,j}^{\mu,\lambda}(v)$, which is wrong.
Therefore $S_{b_m,j}(A_{b_m..j})f_{\mu,\lambda}=0$ for any $m=1,\ldots,N$.

We shall use this fact to prove by downward induction on $u=1,\ldots,N+1$
the following property:
\begin{equation}\label{equation:3}
\begin{array}{l}
\mbox{for any $k=1+\delta_{j=n}(n-1),\ldots,n$, there is a weakly increasing}\\
\mbox{injection $d_k:[b_u..c_u]\cup\cdots\cup[b_N..c_N]\to(i..j)$ such that}\\
\mbox{$B^{\mu,\lambda,k}(x,d_k(x))=0$ for any admissible $x$.}
\end{array}
\end{equation}
This is obviously true for $u=N+1$.
Therefore, we suppose that
$1\le u\le N$ and property~(\ref{equation:3}) is proved for greater $u$.
Fix an arbitrary $k=1+\delta_{j=n}(n-1),\ldots,n$.
Since $S_{b_u,j}(A_{b_u..j})f_{\mu,\lambda}=0$,
the inductive hypothesis asserting that
the current lemma is true for smaller values of $|\bar A|$ implies that
$(b_u..j)\setminus A_{b_u..j}$ satisfies $\pi_{b_u,j}^{\mu,\lambda}(v)$ for some $v$.
As a consequence, there is a weakly increasing injection
$e_k:[b_u..c_u]\cup\cdots\cup[b_{u+w-1}..c_{u+w-1}]\to[b_u..b_{u+w}-1)$
such that $B^{\mu,\lambda,k}(x,d_k(x))=0$ for any admissible $x$
(here $w=v-1+\delta_{b_u=c_u}$ and $b_{N+1}=j+1$).
The inductive hypothesis asserting that property~(\ref{equation:3})
holds for $u+w$ allows us to extend $e_k$ to the required injection $d_k$.
Thus property~(\ref{equation:3}) is proved.

Take any $k=i+\delta_{j=n}(j-i),\ldots,j$.
Applying property~(\ref{equation:3}) for $u=1$,
the fact that $x_k$ follows from $x_j$,
and Proposition~\ref{proposition:2}, we get
$$
0=K^{\mu,\lambda}(x_k)=K_{i,j}^{\mu,\lambda,k}(A)=
\prod_{t\in[i..j)\setminus\im d_k}B^{\mu,\lambda,k}(i,t).
$$
Therefore, there is $t'\in[i..j)\setminus\im d_k$
such that $B^{\mu,\lambda,k}(i,t')=0$.
Putting $\theta_k(t)=d_k(t)$ for
$t\in[b_1..c_1]\cup\cdots\cup[b_N..c_N]$ and
$\theta_k(i)=t'$, we get a map required in Definition~\ref{definition:1}.
This fact together with Remark~\ref{remark:0.5} shows that $\bar A$
satisfies $\pi_{i,j}^{\mu,\lambda}(N+1)$, contrary to assumption.
\endproof

Following~\cite{Kleshchev_gjs11}, we introduce the following sets:
$$
\begin{array}{rcl}
\C^\mu(i,j)      &:=&\{a:i<a<j,C^\mu(i,a)=0\},\\[3pt]
\B^{\mu,\lambda}(i,j)&:=&\{a:i\le a<j,B^{\mu,\lambda}(i,a)=0\},
\end{array}
$$
where $C^\mu(i,a)$ is the residue class of $a-i+\mu_i-\mu_a$ modulo $p$
as in~\cite{Kleshchev_gjs11}.

\begin{theorem}\label{theorem:3}
Let $1\le i<n$.
\begin{enumerate}
\item\label{theorem:3:p:1} Let $A\subset(i..n)$. Then $S_{i,n}(A)f_{\mu,\lambda}$
is a non-zero $U(n-1)$-high weight vector if and only if
there is a weakly increasing injection $d:(i..n)\setminus A\to(i..n)$
such that $B^{\mu,\lambda}(x,d(x))=0$ for any admissible $x$ and
$B^{\mu,\lambda}(i,t)\ne0$ for any $t\in[i..n)\setminus\im d$.
\item\label{theorem:3:p:2} There is some $A\subset(i..n)$ such that
$S_{i,n}(A)f_{\mu,\lambda}$ is a non-zero $U(n-1)$-high weight vector
if and only if
there is a weakly decreasing injection
from $\B^{\mu,\lambda}(i,n)$ to $\C^\mu(i,n)$.
\end{enumerate}
\end{theorem}
\Proof\ref{theorem:3:p:1} It is clear
from~\cite[4.11(ii)]{Brundan_quantum},
Theorem~\ref{theorem:2} and
Proposition~\ref{proposition:2} that $S_{i,n}(A)f_{\mu,\lambda}$
is a non-zero $U(n-1)$-high weight vector for such $A$.
Conversely, if $S_{i,n}(A)f_{\mu,\lambda}$ is a non-zero
$U(n-1)$-high weight vector
then, arguing as in the ``only if part'' of Theorem~\ref{theorem:2},
we get that there is a weakly increasing injection
$d:(i..n)\setminus A\to(i..n)$ such that $B^{\mu,\lambda}(x,d(x))=0$
for any admissible $x$.
Now by Proposition~\ref{proposition:2}, we have
$0\ne K^{\mu,\lambda}(i,n)(A)=\prod_{t\in[i..n)\setminus\im d}B^{\mu,\lambda}(i,t)$.

\ref{theorem:3:p:2} If $\epsilon$ is such an injection,
then it suffices to put $A=(i..n)\setminus\im\epsilon$,
take for $d$ the inverse map of $\epsilon$ and apply part~\ref{theorem:3:p:1}.
Conversely, let $S_{i,n}(A)f_{\mu,\lambda}$ be a non-zero $U(n-1)$-high weight
vector for some $A\subset(i..n)$ and let $d$ be an injection,
whose existence is claimed by part~\ref{theorem:3:p:1}.
Now the result follows from the following
two observations: $\B^{\mu,\lambda}(i,n)\subset\im d$;
$d(x)\in\B^{\mu,\lambda}(i,n)$ implies $x\in\C^\mu(i,n)$.\endproof

\remark\label{remark:1} If we obtain a non-zero $U(n-1)$-high weight vector
in Theorem~\ref{theorem:3}, then it is a scalar multiple of
$f_{\nu,\lambda}$, where $\nu=\mu-\varepsilon_i$ and
$\varepsilon_i=(0^{i-1},1,0^{n-1-i})$.

\section{Moving one node}

\begin{definition}\label{definition:2}
Let $1\le i<j-1<n-1$, $M\subset(i..j-1)$ and
$M=[b_1..c_1]\cup\cdots\cup[b_N..c_N]$ be the decomposition of $M$
into the union of connected components. We say that $M$ satisfies
the condition $\bar\pi_{i,j}^{\mu,\lambda}(v)$ if $1\le v\le N+1$
and for any $k=1,\ldots,j-1$ there exists a weakly increasing
injection
$\theta_k:\{i\}\cup[b_1..c_1]\cup\cdots\cup[b_{v-1}..c_{v-1}]\to[i..b_v-1)$
such that $B^{\mu,\lambda,k}(x,\theta_k(x))=0$ for any admissible
$x$, where we assume $b_{N+1}=j+1$.
\end{definition}

\remark\label{remark:0.75}
If in the above definition for some $k=1,\ldots,j-1$, the
inequality $\theta_k(x)<k$ holds for any admissible $x$, then the
maps
$\theta_l:\{i\}\cup[b_1..c_1]\cup\cdots\cup[b_{v-1}..c_{v-1}]\to[i..b_v-1)$
for $k<l\le n$ such that $B^{\mu,\lambda,l}(x,\theta_l(x))=0$ for
any admissible $x$, can be defined equal to $\theta_k$.

Indeed, it follows from Remark~\ref{remark:0.5} that for $k<l\le
n$ we have
$B^{\mu,\lambda,l}(x,\theta_k(x))=B^{\mu,\lambda,k}(x,\theta_k(x))=0$
for any admissible $x$.
In particular (taking $k=j-1$), we obtain that for $v\le N$ the
set $M$ (that consists of $N$ connected components) satisfies the
condition $\bar\pi_{i,j}^{\mu,\lambda}(v)$ if and only if it
satisfies the condition $\pi_{i,j}^{\mu,\lambda}(v)$.

\begin{theorem}\label{theorem:4}
Let $1\le i<j-1<n-1$ and $A\subset(i..j)$ such that $j-1\in A$.
Then $E_{j-1}S_{i,j}(A)f_{\mu,\lambda}=0$ if and only if
$(i..j-1)\setminus A$ satisfies $\bar\pi_{i,j}^{\mu,\lambda}(v)$
for some $v$.
\end{theorem}
\Proof
Let $\bar A=(i..j)\setminus A$ and
$\bar A=[b_1..c_1]\cup\cdots\cup[b_N..c_N]$ be
the decomposition into the union of connected components.
We put $x_k=\bigl((i,k,j,A)\bigl)$ for brevity.

We prove the theorem by induction on $|\bar A|$. Suppose $\bar
A=\emptyset$. Then all the sequences following from $x_{j-1}$ are
$x_k$, where $i\le k\le j-1$. By Theorem~\ref{theorem:1},
$\Phi(x_{j-1})f_{\mu,\lambda}=0$ if and only if
$K^{\mu,\lambda}(x_k)=0$ for any $k=i,\ldots,j-1$. Applying
Proposition~\ref{proposition:2}, we see that
$\Phi(x_{j-1})f_{\mu,\lambda}=0$ if and only if for any
$k=i,\ldots,j-1$ there is $t_k\in[i..j)$ such that
$B^{\mu,\lambda,k}(i,t_k)=0$. In view of Remark~\ref{remark:0.5},
this assertion is equivalent to $\bar\pi_{i,j}^{\mu,\lambda}(1)$.

Now suppose that $\bar A\ne\emptyset$ and that the theorem holds for
smaller values of $|\bar A|$.

{\it``If part''}. If $v\le N$ then $\bar A$ satisfies
$\pi_{i,j}^{\mu,\lambda}(v)$ by Remark~\ref{remark:0.75}. Hence by
Theorem~\ref{theorem:2}, we have $S_{i,j}(A)f_{\mu,\lambda}=0$ and
the desired result follows.

So we shall consider the case $v=N+1$. For any $m=1,\ldots,N$, the
elements $E_{b_m-1}$ and $E_{j-1}$ commute and
by~\cite[4.11(ii)]{Brundan_quantum} we have $E_{b_m-1}E_{j-1}$
$S_{i,j}(A)f_{\mu,\lambda}= -S_{i,b_m-1}(A_{i..b_m-1})$
$E_{j-1}S_{b_m,j}(A_{b_m..j})f_{\mu,\lambda}$. Note that
\begin{equation}\label{equation:4}
\begin{array}{l}
A_{i..b_m-1}=(i..b_m-1)\setminus\bigl([b_1..c_1]\cup\cdots\cup[b_{m-1}..c_{m-1}]\bigl),\\[3pt]
A_{b_m..j}=(b_m..j)\setminus\bigl((b_m..c_m]\cup\cdots\cup[b_N..c_N]\bigl).
\end{array}
\end{equation}
Obviuosly, the set $(b_m..c_m]\cup\cdots\cup[b_N..c_N]$ satisfies
$\bar\pi_{b_m,j}^{\mu,\lambda}(N+2-m-\delta_{b_m=c_m})$, whence by
the inductive hypothesis
$E_{j-1}S_{b_m,j}(A_{b_m..j})f_{\mu,\lambda}=0$.
Thus we have
\begin{equation}\label{equation:5}
E_{b_m-1}E_{j-1}S_{i,j}(A)f_{\mu,\lambda}=0.
\end{equation}

Let us prove by induction on $s=0,\ldots,j-i-1$ that the conditions
\begin{equation}\label{equation:6}
K_{i,j}^{\mu,\lambda,j-1}(A)=0,\quad \ldots,\quad
K_{i,j}^{\mu,\lambda,j-s}(A)=0,\quad
\Phi(x_{j-1-s})f_{\mu,\lambda}=0
\end{equation}
imply $\Phi(x_{j-1})f_{\mu,\lambda}=0$. It is obviously true for
$s=0$. Suppose that $0<s\le j-i-1$, conditions~(\ref{equation:6})
hold and the assertion is true for smaller values of $s$. By the
inductive hypothesis it suffices to prove that
$\Phi(x_{j-s})f_{\mu,\lambda}=0$. Let
$x_{j-s}\stackrel{l}{\longrightarrow}x'$. We have either
$x'=x_{j-s-1}$ or $l=b_m-1<j-s-1$. Since in the former case
$\Phi(x')f_{\mu,\lambda}=0$ by~(\ref{equation:6}), we shall
consider the latter case. We have
$$
\begin{array}{l}
\Phi(x')f_{\mu,\lambda}=E_{b_m-1}\Phi(x_{j-s})f_{\mu,\lambda}=
E_{b_m-1}E(j-s,j-2)E_{j-1}S_{i,j}(A)f_{\mu,\lambda}\\[6pt]
=E(j-s,j-2)E_{b_m-1}E_{j-1}S_{i,j}(A)f_{\mu,\lambda}=0.
\end{array}
$$
To obtain the last equality, we used~(\ref{equation:5}). Since
$K^{\mu,\lambda}(x_{j-s})=K_{i,j}^{\mu,\lambda,j-s}(A)=0$, we get
$\Phi(x_{j-s})f_{\mu,\lambda}=0$ by Theorem~\ref{theorem:1}.

Note that nothing follows from $x_i$ except itself. Therefore,
applying the above assertion for $s=j-i-1$ and
Theorem~\ref{theorem:1}, we see that to prove
$\Phi(x_{j-1})f_{\mu,\lambda}=0$, it suffices to prove
$K_{i,j}^{\mu,\lambda,k}(A)=0$ for $i\le k\le j-1$. The last
equalities follow from Proposition~\ref{proposition:2}.

{\it``Only if part''}. Suppose $\bar A$ satisfies the condition
$\bar\pi_{i,j}^{\mu,\lambda}(v)$ for no $v$. Multiplying the
equality $\Phi(x_{j-1})f_{\mu,\lambda}=0$ by $E_{b_m-1}$, where
$1\le m\le N$, we get
$S_{i,b_m-1}(A_{i..b_m-1})E_{j-1}S_{b_m,j}(A_{b_m..j})f_{\mu,\lambda}=0$
according to~\cite[4.11(ii)]{Brundan_quantum}. By
Corollary~\ref{corollary:1}, either
$S_{i,b_m-1}(A_{i..b_m-1})f_{\mu,\lambda}=0$ or
$E_{j-1}S_{b_m,j}(A_{b_m..j})f_{\mu,\lambda}=0$. The former case
is impossible since Theorem~\ref{theorem:2} would yield that
$(i..b_m-1)\setminus A_{i..b_m-1}$ satisfies
$\pi_{i,b_m-1}^{\mu,\lambda}(v)$ for some $v\le m$
(see~(\ref{equation:1.25})). But then $\bar A$ would satisfy
$\pi_{i,j}^{\mu,\lambda}(v)$ and thus also would satisfy
$\bar\pi_{i,j}^{\mu,\lambda}(v)$, which is wrong. Therefore
$E_{j-1}S_{b_m,j}(A_{b_m..j})f_{\mu,\lambda}=0$ for any
$m=1,\ldots,N$.

We shall use this fact to prove by downward induction on $u=1,\ldots,N+1$
the following property:
\begin{equation}\label{equation:7}
\begin{array}{l}
\mbox{for any $k=1,\ldots,j-1$, there is a weakly increasing}\\
\mbox{injection $d_k:[b_u..c_u]\cup\cdots\cup[b_N..c_N]\to(i..j)$ such that}\\
\mbox{$B^{\mu,\lambda,k}(x,d_k(x))=0$ for any admissible $x$.}
\end{array}
\end{equation}
This is obviously true for $u=N+1$. Therefore, we suppose that
$1\le u\le N$ and property~(\ref{equation:7}) is proved for
greater $u$. Fix an arbitrary $k=1,\ldots,j-1$. Since
$E_{j-1}S_{b_u,j}(A_{b_u..j})f_{\mu,\lambda}=0$, the inductive
hypothesis asserting that the current lemma is true for smaller
values of $|\bar A|$ implies that $(b_u..j)\setminus A_{b_u..j}$
satisfies $\bar\pi_{b_u,j}^{\mu,\lambda}(v)$ for some $v$. As a
consequence, there is a weakly increasing injection
$e_k:[b_u..c_u]\cup\cdots\cup[b_{u+w-1}..c_{u+w-1}]\to[b_u..b_{u+w}-1)$
such that $B^{\mu,\lambda,k}(x,d_k(x))=0$ for any admissible $x$
(here $w=v-1+\delta_{b_u=c_u}$ and $b_{N+1}=j+1$). The inductive
hypothesis asserting that property~(\ref{equation:7}) holds for
$u+w$ allows us to extend $e_k$ to the required injection $d_k$.
Thus property~(\ref{equation:7}) is proved.

Take any $k=i,\ldots,j-1$.
Applying property~(\ref{equation:7}) for $u=1$,
the fact that $x_k$ follows from $x_j$,
and Proposition~\ref{proposition:2}, we get
$$
0=K^{\mu,\lambda}(x_k)=K_{i,j}^{\mu,\lambda,k}(A)=
\prod_{t\in[i..j)\setminus\im d_k}B^{\mu,\lambda,k}(i,t).
$$
Therefore, there is $t'\in[i..j)\setminus\im d_k$ such that
$B^{\mu,\lambda,k}(i,t')=0$. Putting $\theta_k(t)=d_k(t)$ for
$t\in[b_1..c_1]\cup\cdots\cup[b_N..c_N]$ and $\theta_k(i)=t'$, we
get a map required in Definition~\ref{definition:2}. This fact
together with Remark~\ref{remark:0.5} shows that $\bar A$
satisfies $\bar\pi_{i,j}^{\mu,\lambda}(N+1)$, contrary to
assumption.
\endproof

Following~\cite{Kleshchev_gjs11}, we introduce the following sets:
$$
\begin{array}{rcl}
\C^\mu(i,j)      &:=&\left\{a:i<a<j,a-i+\mu_i-\mu_a\equiv0\pmod p\right\},\\[3pt]
\B^{\mu,\lambda,k}(i,j)&:=&{\{}a:i\le
a<j,B^{\mu,\lambda,k}(i,a)=0{\}}.
\end{array}
$$
We shall abbreviate
$B^\mu(i,a)=B^{\mu,\mu}(i,a)$ and
$\B^\mu(i,j)=\B^{\mu,\mu}(i,j)$.
It follows from Remark~\ref{remark:0.5} that
\begin{equation}\label{equation:8}
\B^{\mu,\lambda,k}(i,j)=
\B^{\mu,\lambda}(i,k) \cup\bigl(\B^\mu(i,j)\cap[k..j)\bigl)
\end{equation}

\begin{theorem}\label{theorem:5}
Let $1\le i<j-1<n-1$.
\begin{enumerate}
\item\label{theorem:5:p:1} Let $A\subset(i..j)$. Then
$S_{i,j}(A)f_{\mu,\lambda}$ is a non-zero $U(n-1)$-high weight
vector if and only if $j-1\in A$, for each $k=1,\ldots,j-1$ there
is a weakly increasing injection $\theta_k:[i..j)\setminus
A\to[i..j)$ such that $B^{\mu,\lambda,k}(x,\theta_k(x))=0$ for any
admissible $x$ and there is a weakly increasing injection
$d:(i..j)\setminus A\to(i..j)$ such that
$B^{\mu,\lambda}(x,d(x))=0$ for any admissible $x$ and
$B^{\mu,\lambda}(i,t)\ne0$ for any $t\in[i..j)\setminus\im d$.
\item\label{theorem:5:p:2} There is some $A\subset(i..j)$ such
that $S_{i,j}(A)f_{\mu,\lambda}$ is a non-zero $U(n-1)$-high
weight vector if and only if there are a weakly decreasing
injection
$\epsilon:\B^{\mu,\lambda}(i,j)\to\C^\mu(i,j-1)$ and weakly
increasing injections
$\theta_k:\{i\}\cup\im\epsilon\to\B^{\mu,\lambda,k}(i,j)$ for any
$k=1,\ldots,j-1$. \item\label{theorem:5:p:3} There is some
$A\subset(i..j)$ such that $S_{i,j}(A)f_{\mu,\lambda}$ is a
non-zero $U(n-1)$-high weight vector if and only if
$j-1\in\B^\mu(i,j)$ {\rm(}i.e. $B^\mu(i,j-1)=0${\rm)},
$j-1\notin\B^{\mu,\lambda}(i,j)$ {\rm(}i.e.
$B^{\mu,\lambda}(i,j-1)\ne0${\rm)}, there are a weakly decreasing
and a weakly increasing injections from $\B^{\mu,\lambda}(i,j-1)$
to $\C^\mu(i,j-1)$ and to $\B^\mu(i,j-1)$ respectively.
\end{enumerate}
\end{theorem}
\Proof\ref{theorem:5:p:1} It is clear
from~\cite[4.11(ii)]{Brundan_quantum}, Theorems~\ref{theorem:2}
and~\ref{theorem:4} and Proposition~\ref{proposition:2} that
$S_{i,j}(A)f_{\mu,\lambda}$ is a non-zero $U(n-1)$-high weight
vector for such $A$.

Conversely, let $S_{i,j}(A)f_{\mu,\lambda}$ be a non-zero
$U(n-1)$-high weight vector. Suppose that $j-1\notin A$. Since
$\bigl((i,j,j,A)\bigl)\stackrel{j-1}{\longrightarrow}\bigl((i,j-1,j,A)\bigl)$,
we have by Theorem~\ref{theorem:1} that
$K_{i,j}^{\mu,\lambda,j}(A)\ne0$ and
$K_{i,j}^{\mu,\lambda,j-1}(A)=0$. However, it is impossible since
by Lemma~\ref{lemma:4}\ref{lemma:4:p:1} and
Remark~\ref{remark:0.5}, we have
$K_{i,j}^{\mu,\lambda,j}(A)=K_{i,j-1}^{\mu,\lambda}(A)=K_{i,j}^{\mu,\lambda,j-1}(A)$.
Thus we have proved that $j-1\in A$.

Arguing as in the ``only if part'' of Theorem~\ref{theorem:2}, we
get that for each $k=1,\ldots,j$ there is a weakly increasing
injection $d_k:(i..j)\setminus A\to(i..j)$ such that
$B^{\mu,\lambda,k}(x,d_k(x))=0$ for any admissible $x$. By
Theorem~\ref{theorem:1}, we have $K_{i,j}^{\mu,\lambda}(A)\ne0$.
Hence by Proposition~\ref{proposition:2}, we have
$B^{\mu,\lambda}(i,t)\ne0$ for any $t\in[i..j)\setminus\im d_k$.
Since each sequence $\bigl((i,k,j,A)\bigl)$, where
$k=i,\ldots,j-1$, follows from $\bigl((i,j,j)\bigl)$ we have
$K_{i,j}^{\mu,\lambda,k}(A)=0$ for each $k=1,\ldots,j-1$. Applying
Proposition~\ref{proposition:2}, we get the required maps
$\theta_1,\ldots,\theta_{j-1}$.

\ref{theorem:5:p:2} If $\epsilon$ and $\theta_1,\ldots,\theta_{j-1}$ are
such injections,
then it suffices to put $A=(i..j)\setminus\im\epsilon$,
take for $d$ the inverse map of $\epsilon$ and apply part~\ref{theorem:5:p:1}.

Conversely, let $S_{i,j}(A)f_{\mu,\lambda}$ be a non-zero
$U(n-1)$-high weight vector for some $A\subset(i..j)$ and let $d$
and $\theta_1,\ldots,\theta_{j-1}$ be injections, whose existence
is claimed by part~\ref{theorem:5:p:1}. Note that the following
two facts: $\B^{\mu,\lambda}(i,j)\subset\im d$;
$d(x)\in\B^{\mu,\lambda}(i,j)$ implies $x\in\C^\mu(i,j)$. Now we
define $\epsilon(d(x)):=x$ for $x\in
d^{-1}(\B^{\mu,\lambda}(i,j))$. Observing that
$\im\epsilon=d^{-1}(\B^{\mu,\lambda}(i,j))\subset(\C^\mu(i,j-1))\cap((i..j)\setminus
A)$ completes the proof.

\ref{theorem:5:p:3} Let $j-1\in\B^\mu(i,j)$,
$j-1\notin\B^{\mu,\lambda}(i,j)$ and
$\epsilon:\B^{\mu,\lambda}(i,j-1)\to\C^\mu(i,j-1)$ and
$\tau:\B^{\mu,\lambda}(i,j-1)\to\B^\mu(i,j-1)$ be a weakly
decreasing and a weakly increasing injections respectively. We
have $\B^{\mu,\lambda}(i,j)=\B^{\mu,\lambda}(i,j-1)$. Thus it
remains to define injections $\theta_1,\ldots,\theta_{j-1}$. For
$x\in\{i\}\cup\im\epsilon$ and $k=1,\ldots,j-1$, we put
$$
\theta_k(x)=\left\{
\begin{array}{ll}
j-1&\mbox{ if }x=i;\\
\epsilon^{-1}(x)&\mbox{ if }i<x\mbox{ and }\epsilon^{-1}(x)<k;\\
\tau(\epsilon^{-1}(x))&\mbox{ if }i<x\mbox{ and }\epsilon^{-1}(x)\ge k;
\end{array}
\right.
$$
One can easily verify with the help of~(\ref{equation:8}) that
$\epsilon,\theta_1,\ldots,\theta_{j-1}$ thus defined satisfy the conditions from
part~\ref{theorem:5:p:2}.

Conversely, let $\epsilon,\theta_1,\ldots,\theta_{j-1}$ be as in
part~\ref{theorem:5:p:2}. For $k=1,\ldots,j-1$, we have
$|\B^{\mu,\lambda,k}(i,j)|\ge|\im\theta_k|=|\{i\}\cup\im\epsilon|=1+|\B^{\mu,\lambda}(i,j)|$.
Taking $k=j-1$ and applying~(\ref{equation:8}), we get
$$
\begin{array}{l}
|\B^{\mu,\lambda}(i,j-1)|+|\B^\mu(i,j)\cap\{j-1\}|=|\B^{\mu,\lambda,j-1}(i,j)|\\[3pt]
\ge1+|\B^{\mu,\lambda}(i,j)|=1+|\B^{\mu,\lambda}(i,j-1)|+|\B^{\mu,\lambda}(i,j)\cap\{j-1\}|.
\end{array}
$$
Hence
$|\B^\mu(i,j)\cap\{j-1\}|=1+|\B^{\mu,\lambda}(i,j)\cap\{j-1\}|$,
whence $j-1\in\B^\mu(i,j)$ and $j-1\notin\B^{\mu,\lambda}(i,j)$.
Next for any $k=1,\ldots,j-1$, we have
$$
\begin{array}{l}
1+|\B^{\mu,\lambda}(i,k)|+|\B^{\mu,\lambda}(i,j-1)\cap[k..j-1)|=
1+|\B^{\mu,\lambda}(i,j)|\le|\B^{\mu,\lambda,k}(i,j)|\\[3pt]
=|\B^{\mu,\lambda}(i,k)|+|\B^\mu(i,j-1)\cap[k..j-1)|+1.
\end{array}
$$
Hence
$|\B^{\mu,\lambda}(i,j-1)\cap[k..j-1)|\le|\B^\mu(i,j-1)\cap[k..j-1)|$
for any $k=1,\ldots,j-1$ and by~\cite[2.2]{Brundan_quantum} there
is a weakly increasing injection
$\tau:\B^{\mu,\lambda}(i,j-1)\to\B^\mu(i,j-1)$. \endproof

\begin{theorem}\label{theorem:6}
Part~\ref{theorem:5:p:3} of Theorem~\ref{theorem:5} remains true for
$1<j=i+1<n$.
\end{theorem}
\Proof
Indeed, $S_{i,i+1}(\emptyset)=F_{i,i+1}$ is a non-zero
$U(n-1)$-high weight vector if and only if
$\mu_i-\lambda_{i+1}\not\equiv0\pmod p$ and
$\mu_i-\mu_{i+1}\equiv0\pmod p$.
Taking into account $\B^{\mu,\lambda}(i,j-1)=\emptyset$,
$B^{\mu,\lambda}(i,j-1)=\mu_i-\lambda_{i+1}+p\Z$ and
$B^\mu(i,j-1)=\mu_i-\mu_{i+1}+p\Z$, we obtain the required result.
\endproof.


\begin{thebibliography}{1}\label{bibl}

\bibitem{Brundan_quantum}
J. Brundan, Modular branching rules and the Mullineux map for Hecke algebras of
  type A, {\em Proc. London Math. Soc.}, {\bf 77} (1998), n. 3, 551--581.

\bibitem{Carter1}
R.W. Carter, Raising and lowering operators for ${\mathfrak sl}_n$, with
  application to orthogonal bases of ${\mathfrak sl}_n$-modules. In Arcata
  conference on representations of finite groups, {\em Proc. Simp. Pure Math.},
  {\bf 47} (1987), 351--366.

\bibitem{Jantzen1}
J.C. Jantzen, Representations of algebraic groups, Pure and Applied
  Mathematics, 131, {\em Academic Press, Inc.}, Boston, MA, 1987.

\bibitem{Kleshchev_gjs11}
A. Kleshchev, J. Brundan and I. Suprunenko, Semisimple restrictions from ${\rm
  GL}(n)$ to ${\rm GL}(n-1)$, {\em J. reine angew. Math.}, {\bf 500} (1998),
  83--112.

\bibitem{Kleshchev2}
A. Kleshchev, Branching rules for modular representations of symmetric
  groups. II, {\em J. Reine Angew. Math.}, {\bf 459} (1995), 163--212.

\end{thebibliography}
\end{document}